\documentclass[aps,prb,superscriptaddress,preprint,longbibliography]{revtex4-2}

\usepackage[utf8]{inputenc}
\usepackage[T1]{fontenc}
\usepackage[english]{babel}

\usepackage{amsmath}
\usepackage{amssymb}
\usepackage{bbold}

\usepackage{graphicx} 

\usepackage[colorlinks=true, urlcolor=blue, linkcolor=blue, citecolor=blue]{hyperref}

\newcommand{\dd}{\mathrm{d}}

\newcommand{\pperp}{{\top}}

\usepackage[normalem]{ulem}

\makeatletter
\patchcmd{\frontmatter@abstract@produce}
  {\vskip200\p@\@plus1fil
   \penalty-200\relax
   \vskip-200\p@\@plus-1fil}
  {}
  {}
  {}
\makeatother

\begin{document}

\setlength{\parskip}{.25\baselineskip}
\setlength{\parindent}{0pt}

\title{Derivative of the truncated singular value and eigen decomposition}

\newcommand{\FUB}{Dahlem Center for Complex Quantum Systems and Institut für Theoretische Physik, Freie Universität Berlin, Arnimallee 14, 14195 Berlin, Germany}

\author{Jan Naumann}
\email{j.naumann@fu-berlin.de}
\affiliation{\FUB}

\date{November 18, 2025}

\begin{abstract}
    \noindent Recently developed applications in the field of machine learning and computational physics rely on automatic differentiation techniques, that require stable and efficient linear algebra gradient computations. This technical note provides a comprehensive and detailed discussion of the derivative of the truncated singular and eigenvalue decomposition. It summarizes previous work and builds on them with an extensive description of how to derive the relevant terms. A main focus is correctly expressing the derivative in terms of the truncated part, despite lacking knowledge of the full decomposition.
\end{abstract}

\maketitle

\tableofcontents

\section{Introduction}

In recent years the usage of automatic differentiation, originally emerged in the context of machine learning, to solve problems in quantum physics lead to the requirement of stable and efficient gradients of linear algebra operations.
The example in mind is the framework of variational optimization of two dimensional tensor networks -- for a comprehensive review of the method see Ref.~\cite{10.21468/SciPostPhysLectNotes.86} and the citations therein.
In this technical note a comprehensive discussion of the derivative of the truncated singular and eigen value decompostition (tSVD/tEVD) is presented.
This work aims to summarize the results already presented in the literature in a self-contained and detailed way. The aim of this work is to formulate the results in a form useful for a technical implementation in the context of automatic differentation frameworks as JAX~\cite{jax2018github}.
The discussion is mainly based on the work from Ref.~\cite{wan2019automaticdifferentiationcomplexvalued, wilde2022scalablylearningquantummanybody, francuz2023stableefficientdifferentiationtensor,aa07_comput_eigen_eigen_deriv_gener}.
Previous works usually assume that one has access to the full singular/eigen decomposition including all singular/eigen values and vectors.
However, as Francuz et al.~\cite{francuz2023stableefficientdifferentiationtensor} pointed out, a correction to the derivative is necessary when only a truncated decomposition is computed and the derivative is obtained from the the kept part.
This note is structured in two main parts, one for the truncated singular value decomposition (tSVD), one for the truncated eigen decomposition (tEVD).
For the tSVD a detailed comparison of the two cases where first a full SVD is calculated and then truncated and second a iterative solver is used to only calculate the kept part is provided.
For the eigen decomposition only the second case is discussed since the other one can be easily derived from it.

\section{Truncated singular value decomposition (SVD)}

\subsection{General definition}

The full singular value decomposition (SVD) of a matrix $A \in \mathbb{C}^{n \times m}$ is expressed as
\begin{align}
    A &= U_f S_f V_f^\dagger , \label{eq:full_svd}
\end{align}
with $S_f \in \mathbb{R}^{r \times r}$ diagonal and positive semi-definite and $U_f \in \mathbb{C}^{n \times r}$ as well as $V_f^\dagger \in \mathbb{C}^{r \times m}$ for $r = \min(n,m)$. $U_f$ and $V_f$ are unitaries if $n = m$, otherwise one matrix is unitary ($U_f$ for $n < m$ or $V_f$ for $n > m$) and the other is semi-unitary with only the property $U_f^\dagger U_f$ or $V_f^\dagger V_f = \mathbb{1}_r$ remaining. Here $\mathbb{1}_r$ denotes the identity matrix of size $r \times r$.

To truncate the singular values down to $t < r$ values, we split the elements into the truncated space and its orthogonal counterpart also known as the discarded space. This can be written as block matrices of the form
\begin{align}
    U_f &= \left(
    \begin{array}{c|c}
        U & U_\perp
    \end{array}
    \right),\\
    S_f &= \left(
    \begin{array}{c|c}
        S & 0\\
        \hline
        0 & S_\perp
    \end{array}
    \right),\\
    V^\dagger_f &= \left(
    \begin{array}{c}
        V^\dagger \\
        \hline
        V^\dagger_\perp
    \end{array}
    \right),
\end{align}
with $U \in \mathbb{C}^{n \times t}$, $U_\perp \in \mathbb{C}^{n \times (r-t)}$, $S \in \mathbb{R}^{t \times t}$, $S_\perp \in \mathbb{R}^{(r-t) \times (r-t)}$, $V^\dagger \in \mathbb{C}^{t \times m}$ and $V^\dagger_\perp \in \mathbb{C}^{(r-t) \times m}$.
Applying these relations to Eq.~\eqref{eq:full_svd} we obtain easily the adapted relation
\begin{align}
    A &= U S V^\dagger + U_\perp S_\perp V_\perp^\dagger . \label{eq:svd_trunc}
\end{align}
From the (orthonormality) properties of the full matrices we can obtain some relations for the block matrices as well,
\begin{align}
    U^\dagger U = V^\dagger V &= \mathbb{1}_t ,  &
    U_\perp^\dagger U_\perp = V_\perp^\dagger V_\perp &= \mathbb{1}_{r-t} , \label{eq:prop_1} \\
    U_\perp^\dagger U &= \mathbb{0}_{(r-t)\times t} , & 
    U^\dagger U_\perp &= \mathbb{0}_{t\times (r-t)} , \label{eq:prop_2} \\
    V_\perp^\dagger V &= \mathbb{0}_{(r-t)\times t} , & 
    V^\dagger V_\perp &= \mathbb{0}_{t\times (r-t)} , \label{eq:prop_3} \\
    U^\dagger U \oplus \mathbb{0}_{r-t} + \mathbb{0}_{t} \oplus U_\perp^\dagger U_\perp &= \mathbb{1}_{r}, & V^\dagger V \oplus \mathbb{0}_{r-t} + \mathbb{0}_{t} \oplus V_\perp^\dagger V_\perp &= \mathbb{1}_{r},
\end{align}
with $\oplus$ indicating the matrix direct sum ($A \oplus B = \operatorname{blockdiag}(A, B)$), $\mathbb{0}_{t}$ the square matrix of size $t \times t$ consisting only of zeros and $\mathbb{0}_{x \times y}$ the matrix of size $x \times y$ consiting only of zeros.
The last line can be compactly written as $\overline{U_\perp^\dagger U_\perp} = \mathbb{1}_r - \overline{U^\dagger U}$ and $\overline{V_\perp^\dagger V_\perp} = \mathbb{1}_r - \overline{V^\dagger V}$, where we extend the matrices to the bigger vector space $\mathbb{C}^{r \times r}$ in the way $\overline{U^\dagger U} \equiv U^\dagger U \oplus \mathbb{0}_{r-t}$ and $\overline{U_\perp^\dagger U_\perp} \equiv \mathbb{0}_{t} \oplus U_\perp^\dagger U_\perp$ (the relations for $V$ and $V_\perp$ are defined analogously).

\subsubsection*{Square case}

In this section we assume that the original matrix $A$ is square and thus $r = m = n$. Then we can use the (unitarity) property, $U_f U_f^\dagger = V_f V_f^\dagger = \mathbb{1}_r$.
Applying this to our split matrices we obtain
\begin{align}
    U U^\dagger + U_\perp U_\perp^\dagger &= \mathbb{1}_r , &
    V V^\dagger + V_\perp V_\perp^\dagger &= \mathbb{1}_r . \label{eq:u_u_dagger_sum}
\end{align}

\subsubsection*{Non-square case}\label{sec:tsvd_non_square}

If our input matrix $A \in \mathbb{C}^{n \times m}$ is not square, then either $U_f$ or $V_f$ is not unitary anymore, but only a semi-unitary.
For the discussion, we assume that $n > m$ and thus $U_f$ is only a semi-unitary as well as $r = m$.
The results for the opposite case $n < m$ can be found analogously.
For the derivative calculation below, it will be necessary to extend $U_f$ to a full unitary in the $\mathbb{C}^{n \times n}$ space.
To do so, we introduce another matrix $U_\pperp \in \mathbb{C}^{n \times (n-m)}$ whose columns are orthonormal vectors orthogonal to the (left-)singular vectors of $U_f$.
A full unitary $U_f' \in \mathbb{C}^{n \times n}$ can be expressed by (horizontally) stacking both matrices $U'_f =\left(\begin{array}{c|c}U_f & U_\pperp\end{array}\right)$.

For the tSVD, this modifies our block matrices for a complete ansatz to the following form,
\begin{align}
    \widetilde{U_f} &= \left(
    \begin{array}{c|c|c}
        U & U_\perp & U_\pperp
    \end{array}
    \right),\\
    \widetilde{S_f} &= \left(
    \begin{array}{c|c|c}
        S & 0 & 0\\
        \hline
        0 & S_\perp & 0\\
        \hline
        0 & 0 & 0
    \end{array}
    \right),\\
    \widetilde{V_f} &= \left(
    \begin{array}{c|c|c}
        V & V_\perp & 0
    \end{array}
    \right).
\end{align}
We want to note that the orthonormal extension matrix $U_\pperp$ obviously fulfills the orthogonality property to the other two parts as well,
\begin{align}
    U_\pperp^\dagger U_\pperp &= \mathbb{1}_{n-m} , &
    U_\pperp^\dagger U &= \mathbb{0}_{(n-m)\times t}, &
    U_\pperp^\dagger U_\perp &= \mathbb{0}_{(n-m)\times (r-t)} . \label{eq:non_sqaure_orth_prop_1}
\end{align}

Since we have two unitaries again, we can find two completeness relations similar to the square case,
\begin{align}
    U U^\dagger + U_\perp U_\perp^\dagger + U_\pperp U_\pperp^\dagger  &= \mathbb{1}_n , &
    V V^\dagger + V_\perp V_\perp^\dagger &= \mathbb{1}_m . \label{eq:non_square_u_u_dagger_sum}
\end{align}
This relation also allows us to easily express $U_\pperp U_\pperp^\dagger$ in terms of known matrices,
\begin{align}
    U_\pperp U_\pperp^\dagger  &= \mathbb{1}_n - U U^\dagger - U_\perp U_\perp^\dagger. \label{eq:U_pperp_U_pperp_dagger}
\end{align}

\subsection{General derivative}

From Eq.~\eqref{eq:svd_trunc} we can obtain the differential of the SVD relation for the truncated case,
\begin{align}
    \dd A &= \dd U S V^\dagger + U \dd S V^\dagger + U S \dd V^\dagger + \dd U_\perp S_\perp V_\perp^\dagger + U_\perp \dd S_\perp V_\perp^\dagger + U_\perp S_\perp \dd V_\perp^\dagger. \label{eq:diff_svd_trunc}
\end{align}

The properties from Eqs.~(\ref{eq:prop_1},~\ref{eq:prop_2},~\ref{eq:prop_3}) can be used to find the following orthogonality-preserving differential conditions relations
\begin{align}
    \dd U^\dagger U + U^\dagger \dd U &= 0 , &
    \dd V^\dagger V + V^\dagger \dd V &= 0 , \label{eq:diff_prop_1} \\
    \dd U_\perp^\dagger U + U_\perp^\dagger \dd U &= 0 , &
    \dd U^\dagger U_\perp + U^\dagger \dd U_\perp &= 0 , \label{eq:diff_prop_2} \\
    \dd V_\perp^\dagger V + V_\perp^\dagger \dd V &= 0 , &
    \dd V^\dagger V_\perp + V^\dagger \dd V_\perp &= 0 . \label{eq:diff_prop_3}
\end{align}

In the discussion, we distinguish between the square and the non-square case for the sake of clarity, since the latter requires a more thorough analysis.

\subsection{Derivative for square case}\label{sec:svd_diff_square}

Since we are interested in how the differentials $\dd U$, $\dd S$ and $\dd V$ are modified by the truncation, we take a look at the $\dd A$ terms split in the both spaces by using Eq.~\eqref{eq:u_u_dagger_sum},
\begin{align}
    \dd A &= U U^\dagger \dd A + (\mathbb{1}_r - U U^\dagger) \dd A = U U^\dagger \dd A + U_\perp U_\perp^\dagger \dd A , \label{eq:split_dA_U} \\
    \dd A &= \dd A V V^\dagger + \dd A (\mathbb{1}_r - V V^\dagger) = \dd A V V^\dagger + \dd A V_\perp V_\perp^\dagger . \label{eq:split_dA_V}
\end{align}
Notice that we are just projecting with $UU^\dagger$ $\left( VV^\dagger \right)$ onto the truncated space spanned by the left (right) singular basis vectors, respectively. This also involves using their orthogonal counterpart, $U_\perp {U_\perp}^{\dagger}$ $\left( V_\perp {V_\perp}^{\dagger}\right)$, which instead projects to the discarded space spanned by the left (right) singular basis vectors, respectively. Combining both previous equations yields
\begin{subequations}
  \begin{align}
    \dd A &= U U^\dagger \dd A V V^\dagger \label{eq:svd_dA_projected_1}\\
          &+ U_\perp U_\perp^\dagger \dd A V_\perp V_\perp^\dagger \label{eq:svd_dA_projected_2}\\
          &+ U U^\dagger \dd A V_\perp V_\perp^\dagger \label{eq:svd_dA_projected_3}\\
          &+ U_\perp U_\perp^\dagger \dd A V V^\dagger . \label{eq:svd_dA_projected_4}
  \end{align}
\end{subequations}
The first two lines \eqref{eq:svd_dA_projected_1} and \eqref{eq:svd_dA_projected_2} correspond to the parts of the differential that only live in the truncated and the discarded space, respectively. The last two lines \eqref{eq:svd_dA_projected_3} and \eqref{eq:svd_dA_projected_4} contain the cross terms. We can now insert Eq.~\eqref{eq:diff_svd_trunc} into Eqs.~(\ref{eq:svd_dA_projected_1}-\ref{eq:svd_dA_projected_4}), whereby we use the orthogonality properties from Eqs.~(\ref{eq:prop_2},~\ref{eq:prop_3}),
\begin{subequations}
  \begin{align}
    \dd A &= U \dd S V^\dagger + U U^\dagger \dd U S V^\dagger + U S \dd V^\dagger V V^\dagger \\
          &+ U_\perp \dd S_\perp V_\perp^\dagger + U_\perp U_\perp^\dagger \dd U_\perp S_\perp V_\perp^\dagger + U_\perp S_\perp \dd V_\perp^\dagger V_\perp V_\perp^\dagger \\
          &+ U S \dd V^\dagger V_\perp V_\perp^\dagger + U U^\dagger \dd U_\perp S_\perp V_\perp^\dagger \\
          &+ U_\perp U_\perp^\dagger \dd U S V^\dagger + U_\perp S_\perp \dd V_\perp^\dagger V V^\dagger .
  \end{align}
\end{subequations}
Applying the orthogonality-preserving differential relations from Eqs.~(\ref{eq:diff_prop_2},~\ref{eq:diff_prop_3}), we obtain
\begin{subequations}
  \begin{align}
    \dd A &= U \dd S V^\dagger + U U^\dagger \dd U S V^\dagger + U S \dd V^\dagger V V^\dagger \\
          &+ U_\perp \dd S_\perp V_\perp^\dagger + U_\perp U_\perp^\dagger \dd U_\perp S_\perp V_\perp^\dagger + U_\perp S_\perp \dd V_\perp^\dagger V_\perp V_\perp^\dagger \\
          &+ U S \dd V^\dagger V_\perp V_\perp^\dagger - U \dd U^\dagger U_\perp S_\perp V_\perp^\dagger \\
          &+ U_\perp U_\perp^\dagger \dd U S V^\dagger - U_\perp S_\perp V_\perp^\dagger \dd V V^\dagger .
  \end{align}
  \label{eq:svd_da_not_compact}
\end{subequations}
For compactness, we define the following notation
\begin{align}
    \widetilde{\dd U}_1 &= U^\dagger \dd U , & \widetilde{\dd U}_2 &= U_\perp^\dagger \dd U , \label{eq:compact_U_tilde}\\
    \widetilde{\dd V}_1 &= V^\dagger \dd V , & \widetilde{\dd V}_2 &= V_\perp^\dagger \dd V . \label{eq:compact_V_tilde}
\end{align}
This notation is very convenient, as it allows us to split $\dd U$/$\dd V$ in the same way as we did previously for $\dd A$ in Eqs.~(\ref{eq:split_dA_U},~\ref{eq:split_dA_V}),
\begin{align}
    \dd U &= U U^\dagger \dd U + U_\perp U_\perp^\dagger \dd U = U \widetilde{\dd U}_1 + U_\perp \widetilde{\dd U}_2 , \label{eq:full_d_U} \\
    \dd V &= V V^\dagger \dd V + V_\perp V_\perp^\dagger \dd V = V \widetilde{\dd V}_1 + V_\perp \widetilde{\dd V}_2 . \label{eq:full_d_V}
\end{align}
Rewriting our last expression for $\dd A$, Eq.~\eqref{eq:svd_da_not_compact}, in terms of the compact notation, we obtain
\begin{subequations}
  \begin{align}
    \dd A &= U \dd S V^\dagger + U \widetilde{\dd U}_1 S V^\dagger + U S \widetilde{\dd V}_1^\dagger V^\dagger \\
          &+ U_\perp \dd S_\perp V_\perp^\dagger + U_\perp U_\perp^\dagger \dd U_\perp S_\perp V_\perp^\dagger + U_\perp S_\perp \dd V_\perp^\dagger V_\perp V_\perp^\dagger \\
          &+ U S \widetilde{\dd V}_2^\dagger V_\perp^\dagger - U \widetilde{\dd U}_2^\dagger S_\perp V_\perp^\dagger \label{eq:svd_da_compact_3}\\
          &+ U_\perp \widetilde{\dd U}_2 S V^\dagger - U_\perp S_\perp \widetilde{\dd V}_2 V^\dagger . \label{eq:svd_da_compact_4}
  \end{align}
\end{subequations}

Our goal is now to determine $\dd S$, $\dd U$ and $\dd V$. As shown in Eqs.~(\ref{eq:full_d_U},~\ref{eq:full_d_V}), the latter two Eqs.~(\ref{eq:svd_da_compact_3},~\ref{eq:svd_da_compact_4}) can be described by determining $\widetilde{\dd U}_1$/$\widetilde{\dd U}_2$ and $\widetilde{\dd V}_1$/$\widetilde{\dd V}_2$.
To further simplify the previous expression for $\dd A$, it is convenient to compute the joint action of the (semi-)unitaries from left and right over $\dd A$ further simplifying with the orthogonality properties, which yields
\begin{align}
    U^\dagger \dd A V &= \dd S + \widetilde{\dd U}_1 S + S \widetilde{\dd V}_1^\dagger , \label{eq:U_dagger_dA_V} \\
    U_\perp^\dagger \dd A V_\perp &= \dd S_\perp + U_\perp^\dagger \dd U_\perp S_\perp + S_\perp \dd V_\perp^\dagger V_\perp , \label{eq:svd_dA_split_up_2}\\
    U^\dagger \dd A V_\perp &= S \widetilde{\dd V}_2^\dagger - \widetilde{\dd U}_2^\dagger S_\perp ,  \label{eq:U_dagger_dA_V_perp} \\
    U_\perp^\dagger \dd A V &= \widetilde{\dd U}_2 S - S_\perp \widetilde{\dd V}_2 . \label{eq:U_perp_dagger_dA_V}
\end{align}
The second equation~\eqref{eq:svd_dA_split_up_2} only contains terms living supported in the discarded space, and thus are neglected in the following. 

\subsubsection*{Calculation of $\dd S$, $\widetilde{\dd U}_1$, $\widetilde{\dd V}_1$}

To determine the relation for $\dd S$, we want to note that from Eq.~\eqref{eq:diff_prop_1} follows that $\widetilde{\dd U}_1$ and $\widetilde{\dd V}_1$ are anti-Hermitian, and thus have a diagonal with purely imaginary entries. This property can be used to calculate $\dd S$ by adding Eq.~\eqref{eq:U_dagger_dA_V} and its Hermitian conjugate (note that $\dd S$ has to be real and diagonal), which yields
\begin{align}
    U^\dagger \dd A V + V^\dagger \dd A^\dagger U &= 2 \dd S + \widetilde{\dd U}_1 S - S \widetilde{\dd U}_1 + S \widetilde{\dd V}_1^\dagger - \widetilde{\dd V}_1^\dagger S.
\end{align}
Restricting ourselves solely to the diagonal entries, we can fully determine dS as follows
\begin{align}
     \dd S &= \mathbb{1}_t \circ \frac{1}{2} \left(U^\dagger \dd A V + V^\dagger \dd A^\dagger U\right) ,
\end{align}
where $\circ$ denotes the Hadamard or element-wise  product of matrices, filtering out only the diagonal entries.
Considering now the off-diagonal parts $U^\dagger \mathrm{d} A V$ leads to an expression for the calculation of the off-diagonal elements of $\widetilde{\dd U}_1$ and $\widetilde{\dd V}_1$,
\begin{align}
    \left(U^\dagger \dd A V + V^\dagger \dd A^\dagger U\right)_{\substack{i,j \\ i \neq j}} &= (\widetilde{\dd U}_{1})_{i,j} S_{j,j} - S_{i,i} (\widetilde{\dd U}_{1})_{i,j} + S_{i,i} (\widetilde{\dd V}_{1}^\dagger)_{i,j} - (\widetilde{\dd V}_{1}^\dagger)_{i,j} S_{j,j}. \label{eq:svd_da_conj_offdiagional}
\end{align}
By appropriately multiplying by $S$ on Eq.~\eqref{eq:svd_da_conj_offdiagional} (for clearness marked by a bold $\boldsymbol{S}$ and square symbol~$\boldsymbol{^2}$), we can obtain the following relations for the off-diagonal elements of $\widetilde{\dd U}_1$ and $\widetilde{\dd V}_1$,
\begin{align}
    \left(U^\dagger \dd A V \boldsymbol{S} + \boldsymbol{S} V^\dagger \dd A^\dagger U\right)_{\substack{i,j \\ i \neq j}} &= (\widetilde{\dd U}_1)_{i,j} S^{\boldsymbol{2}}_{j,j} - S^{\boldsymbol{2}}_{i,i} (\widetilde{\dd U}_1)_{i,j} , \\
    \left(\boldsymbol{S} U^\dagger \dd A V + V^\dagger \dd A^\dagger U \boldsymbol{S}\right)_{\substack{i,j \\ i \neq j}} &= S^{\boldsymbol{2}}_{i,i} (\widetilde{\dd V}_1^\dagger)_{i,j} - (\widetilde{\dd V}_1^\dagger)_{i,j} S^{\boldsymbol{2}}_{j,j} ,
\end{align}
where we have used the anti-Hermiticity of $\widetilde{\dd V}_1$, i.e., $\widetilde{\dd V}_1^\dagger = -\widetilde{\dd V}_1$.
Assuming now that we have no degenerate singular values, we can define a matrix $F \in \mathbb{R}^{t \times t}$ with $F_{i,j} = 1/(S^2_{j,j} - S^2_{i,i})$ for $i \neq j$ and zero on the diagonal. Using this matrix, we can rewrite the result for the off-diagonal entries of our differentials as 
\begin{align}
    \bar{\mathbb{1}}_t \circ \widetilde{\dd U}_1 &= F \circ \left(U^\dagger \dd A V S + S V^\dagger \dd A^\dagger U\right) , \\
    \bar{\mathbb{1}}_t \circ \widetilde{\dd V}_1 &= F \circ \left(S U^\dagger \dd A V + V^\dagger \dd A^\dagger U S\right) ,
\end{align}
where $\bar{\mathbb{1}}_t$ denotes the square matrix of size $t \times t$ containing ones in all off-diagonal entries, and zeros along the diagonal.

Looking at the anti-Hermitian part of Eq.~\eqref{eq:U_dagger_dA_V}; i.e., subtracting its corresponding hermitian conjugate version from itself, yields
\begin{align}
    U^\dagger \dd A V - V^\dagger \dd A^\dagger U &= \widetilde{\dd U}_1 S + S \widetilde{\dd U}_1 + S \widetilde{\dd V}_1^\dagger + \widetilde{\dd V}_1^\dagger S.
\end{align}
Now we can solve for the diagonal part of $\widetilde{\dd U}_1 + \widetilde{\dd V}_1$, obtaining
\begin{align}
    \Rightarrow \mathbb{1}_t \circ \left(\widetilde{\dd U}_1 + \widetilde{\dd V}_1^\dagger\right) &= S^{-1} \circ \frac{1}{2} \left(U^\dagger \dd A V - V^\dagger \dd A^\dagger U\right) .
\end{align}
As discussed in Ref.~\cite{wilde2022scalablylearningquantummanybody, francuz2023stableefficientdifferentiationtensor}, any solution of this sum is enough to reconstruct $\dd A$, so we can just add the term to either of the differentials or split it between both. Wilde et al.~\cite{wilde2022scalablylearningquantummanybody} suggested that splitting the term in half yield to smallest rounding error in their numerical analysis.

\subsubsection*{Calculation of $\widetilde{\dd U}_2$, $\widetilde{\dd V}_2$}

A similar approach as for the calculation of the off-diagonal elements above can be used to obtain $\widetilde{\dd U}_2$ and $\widetilde{\dd V}_2$. Let us take a look at the sum of the hermitian conjugate of Eq.~\eqref{eq:U_dagger_dA_V_perp} and Eq.~\eqref{eq:U_perp_dagger_dA_V},
\begin{align}
    U_\perp^\dagger \dd A V + V_\perp^\dagger \dd A^\dagger U &= \widetilde{\dd U}_2 S - S_\perp \widetilde{\dd V}_2 + \widetilde{\dd V}_2 S - S_\perp \widetilde{\dd U}_2 .
\end{align}
Again, we can modify the equation by multiply $S$ and $S_\perp$ from the left and the right as indicted by bold symbols,
\begin{align}
    U_\perp^\dagger \dd A V \boldsymbol{S} + \boldsymbol{S_\perp} V_\perp^\dagger \dd A^\dagger U &= \widetilde{\dd U}_2 S^{\boldsymbol{2}} - S^{\boldsymbol{2}}_\perp \widetilde{\dd U}_2 , \\
    \boldsymbol{S_\perp} U_\perp^\dagger \dd A V + V_\perp^\dagger \dd A^\dagger U \boldsymbol{S} &= \widetilde{\dd V}_2 S^{\boldsymbol{2}} - S^{\boldsymbol{2}}_\perp \widetilde{\dd V}_2.
\end{align}
Defining a matrix $G \in \mathbb{R}^{(r-t) \times t}$ as $G_{i,j} = 1/(S^2_{j, j} - S^2_{\perp_{i, i}})$, and zero elsewhere, both previous equations can be rewritten as
\begin{align}
    \widetilde{\dd U}_2 &= G \circ \left(U_\perp^\dagger \dd A V S + S_\perp V_\perp^\dagger \dd A^\dagger U\right) , \\
    \widetilde{\dd V}_2 &= G \circ \left(S_\perp U_\perp^\dagger \dd A V + V_\perp^\dagger \dd A^\dagger U S\right) .
\end{align}

\subsubsection*{Summary of the results for square case}

To summarize the results for the square case, we have 
\begin{align}
    \dd S &= \mathbb{1}_t \circ \frac{1}{2} \left(U^\dagger \dd A V + V^\dagger \dd A^\dagger U\right) , \\
    \dd U &= U \widetilde{\dd U}_1 + U_\perp \widetilde{\dd U}_2 , \\
    \dd V &= V \widetilde{\dd V}_1 + V_\perp \widetilde{\dd V}_2 ,
\end{align}
with the terms
\begin{align}
    \widetilde{\dd U}_1 &= F \circ \left(U^\dagger \dd A V S + S V^\dagger \dd A^\dagger U\right) + \frac{1}{2} S^{-1} \circ \frac{1}{2} \left(U^\dagger \dd A V - V^\dagger \dd A^\dagger U\right) , \\
    \widetilde{\dd V}_1 &= F \circ \left(S U^\dagger \dd A V + V^\dagger \dd A^\dagger U S\right) + \frac{1}{2} S^{-1} \circ \frac{1}{2} \left(U^\dagger \dd A V - V^\dagger \dd A^\dagger U\right) , \\
    \widetilde{\dd U}_2 &= G \circ \left(U_\perp^\dagger \dd A V S + S_\perp V_\perp^\dagger \dd A^\dagger U\right) , \\
    \widetilde{\dd V}_2 &= G \circ \left(S_\perp U_\perp^\dagger \dd A V + V_\perp^\dagger \dd A^\dagger U S\right) .
\end{align}

\subsection{Derivative for non-square case}

As discussed in Sec.~\ref{sec:tsvd_non_square} we limit our discussion to the case $n > m$. The main modification compared to the square case is that the relation $U U^\dagger + U_\perp U_\perp^\dagger = \mathbb{1}_n$ does not hold anymore, therefore has to be replaced by the new relations described in Eq.~\eqref{eq:non_square_u_u_dagger_sum}.
Thus the ansatz for $\dd U$ from Eq.~\eqref{eq:full_d_U} and the equation for $\dd A$ from Eqs.~(\ref{eq:split_dA_U},~\ref{eq:split_dA_V}) have to be accordingly modified, leading to a new set of equations,
\begin{align}
    \dd A &= U U^\dagger \dd A + (\mathbb{1}_n - U U^\dagger) \dd A = U U^\dagger \dd A + U_\perp U_\perp^\dagger \dd A + U_\pperp U_\pperp^\dagger \dd A , \label{eq:svd_non_sqaure_ansatz_1} \\
    \dd A &= \dd A V V^\dagger + \dd A (\mathbb{1}_m - V V^\dagger) = \dd A V V^\dagger + \dd A V_\perp V_\perp^\dagger , \label{eq:svd_non_sqaure_ansatz_2}\\
    \dd U &= U U^\dagger \dd U + (\mathbb{1}_n - U U^\dagger) \dd U = U U^\dagger \dd U + U_\perp U_\perp^\dagger \dd U + U_\pperp U_\pperp^\dagger \dd U, \label{eq:full_dU_non_square} \\
    \dd V &= V V^\dagger \dd V + (\mathbb{1}_m - V V^\dagger) \dd V = V V^\dagger \dd V + V_\perp V_\perp^\dagger \dd V .
\end{align}
Combining the first two equations~(\ref{eq:svd_non_sqaure_ansatz_1},~\ref{eq:svd_non_sqaure_ansatz_2}) again, as in the square case, yields
\begin{subequations}
  \begin{align}
    \dd A &= U U^\dagger \dd A V V^\dagger \\
          &+ U_\perp U_\perp^\dagger \dd A V_\perp V_\perp^\dagger \\
          &+ U U^\dagger \dd A V_\perp V_\perp^\dagger \\
          &+ U_\perp U_\perp^\dagger \dd A V V^\dagger \\
          &+ U_\pperp U_\pperp^\dagger \dd A V V^\dagger \\
          &+ U_\pperp U_\pperp^\dagger \dd A V_\perp V_\perp^\dagger .
  \end{align}
\end{subequations}
Inserting Eq.~\eqref{eq:diff_svd_trunc} into our relation and using the orthogonality properties from Eqs.~(\ref{eq:prop_1}-\ref{eq:prop_3},~\ref{eq:non_sqaure_orth_prop_1}) we obtain
\begin{subequations}
  \begin{align}
    \dd A &= U \dd S V^\dagger + U U^\dagger \dd U S V^\dagger + U S \dd V^\dagger V V^\dagger \\
          &+ U_\perp \dd S_\perp V_\perp^\dagger + U_\perp U_\perp^\dagger \dd U_\perp S_\perp V_\perp^\dagger + U_\perp S_\perp \dd V_\perp^\dagger V_\perp V_\perp^\dagger \\
          &+ U S \dd V^\dagger V_\perp V_\perp^\dagger + U U^\dagger \dd U_\perp S_\perp V_\perp^\dagger \\
          &+ U_\perp U_\perp^\dagger \dd U S V^\dagger + U_\perp S_\perp \dd V_\perp^\dagger V V^\dagger \\
          &+ U_\pperp U_\pperp^\dagger \dd U S V^\dagger \\
          &+ U_\pperp U_\pperp^\dagger \dd U_\perp S_\perp V_\perp^\dagger.
  \end{align}
\end{subequations}
Using Eqs.~(\ref{eq:diff_prop_2},~\ref{eq:diff_prop_3}), which are still valid, and the compact notation introduced in Eqs.~(\ref{eq:compact_U_tilde},~\ref{eq:compact_V_tilde}) we get the equation
\begin{subequations}
  \begin{align}
    \dd A &= U \dd S V^\dagger + U \widetilde{\dd U}_1 S V^\dagger + U S \widetilde{\dd V}_1^\dagger V^\dagger \\
          &+ U_\perp \dd S_\perp V_\perp^\dagger + U_\perp U_\perp^\dagger \dd U_\perp S_\perp V_\perp^\dagger + U_\perp S_\perp \dd V_\perp^\dagger V_\perp V_\perp^\dagger \\
          &+ U S \widetilde{\dd V}_2^\dagger V_\perp^\dagger - U \widetilde{\dd U}_2^\dagger S_\perp V_\perp^\dagger \\
          &+ U_\perp \widetilde{\dd U}_2 S V^\dagger - U_\perp S_\perp \widetilde{\dd V}_2 V^\dagger \\
          &+ U_\pperp U_\pperp^\dagger \dd U S V^\dagger \\
          &+ U_\pperp U_\pperp^\dagger \dd U_\perp S_\perp V_\perp^\dagger.
  \end{align}
\end{subequations}
This equation is very similar to the ansatz found for the square case except for the last two lines. Therefore, it is a useful idea to again observe the separate components generated by jointly multiplying by the (semi-)unitaries from left and right,
\begin{align}
    U^\dagger \dd A V &= \dd S + \widetilde{\dd U}_1 S + S \widetilde{\dd V}_1^\dagger , \\
    U_\perp^\dagger \dd A V_\perp &= \dd S_\perp + U_\perp^\dagger \dd U_\perp S_\perp + S_\perp \dd V_\perp^\dagger V_\perp , \\
    U^\dagger \dd A V_\perp &= S \widetilde{\dd V}_2^\dagger - \widetilde{\dd U}_2^\dagger S_\perp , \\
    U_\perp^\dagger \dd A V &= \widetilde{\dd U}_2 S - S_\perp \widetilde{\dd V}_2 , \\
    U_\pperp^\dagger \dd A V &= U_\pperp^\dagger \dd U S , \label{eq:U_pperp_dA_V} \\
    U_\pperp^\dagger \dd A V_\perp &= U_\pperp^\dagger \dd U_\perp S_\perp .
\end{align}
Since the first four equations are still unmodified, the procedure to obtain $\dd S$, $\widetilde{\dd U}_1$/$\widetilde{\dd U}_2$ and $\widetilde{\dd V}_1$/$\widetilde{\dd V}_2$ is the same as above. 
It now just remains to address the additional term $U_\pperp U_\pperp^\dagger \dd U$ in Eq.~\eqref{eq:full_dU_non_square}, for which we can look at Eq.~\eqref{eq:U_pperp_dA_V} and rewrite it,
\begin{align}
    U_\pperp U_\pperp^\dagger \dd U &= U_\pperp U_\pperp^\dagger \dd A V S^{-1} .
\end{align}
To avoid calculating $U_\pperp$ explicitly we use Eq.~\eqref{eq:U_pperp_U_pperp_dagger} in the above obtained expression,
\begin{align}
    U_\pperp U_\pperp^\dagger \dd U &= \left(\mathbb{1}_n - U U^\dagger - U_\perp U_\perp^\dagger\right) \dd A V S^{-1} .
\end{align}

For the case $n < m$, the important modifications compared to the square case are 
\begin{align}
    \dd A &= U U^\dagger \dd A + (\mathbb{1}_n - U U^\dagger) \dd A = U U^\dagger \dd A + U_\perp U_\perp^\dagger \dd A , \\
    \dd A &= \dd A V V^\dagger + \dd A (\mathbb{1}_m - V V^\dagger) = \dd A V V^\dagger + \dd A V_\perp V_\perp^\dagger + \dd A V_\pperp V_\pperp^\dagger , \\
    \dd U &= U U^\dagger \dd U + (\mathbb{1}_n - U U^\dagger) \dd U = U U^\dagger \dd U + U_\perp U_\perp^\dagger \dd U, \\
    \dd V &= V V^\dagger \dd V + (\mathbb{1}_m - V V^\dagger) \dd V = V V^\dagger \dd V + V_\perp V_\perp^\dagger \dd V + V_\pperp V_\pperp^\dagger \dd V.
\end{align}
Using the same procedure as above we get the two additional separate components
\begin{align}
    U^\dagger \dd A V_\pperp &= S \dd V^\dagger V_\pperp , \\
    U_\perp^\dagger \dd A V_\pperp &= S_\perp \dd V_\perp^\dagger V_\pperp .
\end{align}
This directly leads to the additional component we need for $\dd V$, which we again analogously rewrite to avoid calculating $V_{\pperp}$,
\begin{align}
    V_\pperp V_\pperp^\dagger \dd V &= V_\pperp V_\pperp^\dagger \dd A^\dagger U S^{-1} = \left(\mathbb{1}_n - V V^\dagger - V_\perp V_\perp^\dagger\right) \dd A^\dagger U S^{-1} .
\end{align}

\subsubsection*{Summary of the results for non-square case}

To summarize the results for the non-square case with $n > m$, we have 
\begin{align}
    \dd S &= \mathbb{1}_t \circ \frac{1}{2} \left(U^\dagger \dd A V + V^\dagger \dd A^\dagger U\right) , \\
    \dd U &= U \widetilde{\dd U}_1 + U_\perp \widetilde{\dd U}_2 + \left(\mathbb{1}_n - U U^\dagger - U_\perp U_\perp^\dagger\right) \dd A V S^{-1}, \\
    \dd V &= V \widetilde{\dd V}_1 + V_\perp \widetilde{\dd V}_2 ,
\end{align}
with the terms
\begin{align}
    \widetilde{\dd U}_1 &= F \circ \left(U^\dagger \dd A V S + S V^\dagger \dd A^\dagger U\right) + S^{-1} \circ \frac{1}{2} \left(U^\dagger \dd A V - V^\dagger \dd A^\dagger U\right) , \\
    \widetilde{\dd V}_1 &= F \circ \left(S U^\dagger \dd A V + V^\dagger \dd A^\dagger U S\right) , \\
    \widetilde{\dd U}_2 &= G \circ \left(U_\perp^\dagger \dd A V S + S_\perp V_\perp^\dagger \dd A^\dagger U\right) , \\
    \widetilde{\dd V}_2 &= G \circ \left(S_\perp U_\perp^\dagger \dd A V + V_\perp^\dagger \dd A^\dagger U S\right) .
\end{align}

For the case $n < m$ we obtain the analog result, 
\begin{align}
    \dd S &= \mathbb{1}_t \circ \frac{1}{2} \left(U^\dagger \dd A V + V^\dagger \dd A^\dagger U\right) , \\
    \dd U &= U \widetilde{\dd U}_1 + U_\perp \widetilde{\dd U}_2 , \\
    \dd V &= V \widetilde{\dd V}_1 + V_\perp \widetilde{\dd V}_2 + \left(\mathbb{1}_n - V V^\dagger - V_\perp V_\perp^\dagger\right) \dd A^\dagger U S^{-1} ,
\end{align}
with the same terms as above.

\subsection{Derivative for iterative methods}

In the discussion so far we assumed that we know the full singular value decomposition and and then truncate it to a smaller subspace. But in practice, there are several methods to calculate only the $t$ usually largest singular values by iterative/Krylov methods, for example, using the Golub-Kahan-Lanczos method~\cite{doi:10.1137/0702016,doi:10.1137/S1064827597327309,doi:10.1137/04060593X}.
As discussed in Ref.~\cite{francuz2023stableefficientdifferentiationtensor} it is still possible to the differentials resulting from the parts in the discarded space.

As above, we can express the singular value decomposition of a matrix $A \in \mathbb{C}^{n \times m}$ as
\begin{align}
    A &= U S V^\dagger + U_\perp S_\perp V_\perp^\dagger. \label{eq:A_tSVD_gkl}
\end{align}

Here we assume that $t < \min(n, m)$. If $t$ is exactly equal to the minimum of $n$ and $m$, we refer to the previously discussed case. We want to note again that the matrices $U$ and $V$ are semi-unitaries fulfilling the property $U^\dagger U = V^\dagger V = \mathbb{1}_t$. The problem in the discussion of this section is that we do not want to calculate $U_\perp$, $S_\perp$ and $V_\perp^\dagger$ explicitly.

The differential of this equation can be still written as expressed in Eq.~\eqref{eq:diff_svd_trunc}. Also, the properties in Eqs.~(\ref{eq:prop_1}-\ref{eq:prop_3},~\ref{eq:diff_prop_1}-\ref{eq:diff_prop_3}) are still true. To start calculating the differentials we have to modify the ansatz to split $\dd A$, $\dd U$ and $\dd V$ up,
\begin{align}
    \dd A &= U U^\dagger \dd A + (\mathbb{1}_n - U U^\dagger) \dd A ,\\
    \dd A &= \dd A V V^\dagger + \dd A (\mathbb{1}_m - V V^\dagger) , \\
    \dd U &= U U^\dagger \dd U + (\mathbb{1}_n - U U^\dagger) \dd U , \label{eq:dU_full_gkl} \\
    \dd V &= V V^\dagger \dd V + (\mathbb{1}_m - V V^\dagger) \dd V . \label{eq:dV_full_gkl}
\end{align}
Combining the first two equations we obtain
\begin{subequations}
  \begin{align}
    \dd A &= U U^\dagger \dd A V V^\dagger \\
          &+ (\mathbb{1}_n - U U^\dagger) \dd A (\mathbb{1}_m - V V^\dagger) \\
          &+ U U^\dagger \dd A V_\perp (\mathbb{1}_m - V V^\dagger) \\
          &+ (\mathbb{1}_n - U U^\dagger) \dd A V V^\dagger .
  \end{align}
\end{subequations}
Inserting Eq.~\eqref{eq:diff_svd_trunc} into our relation and using the orthogonality properties from Eqs.~(\ref{eq:prop_1}-\ref{eq:prop_3},~\ref{eq:non_sqaure_orth_prop_1}), we obtain
\begin{subequations}
  \begin{align}
    \dd A &= U \dd S V^\dagger + U U^\dagger \dd U S V^\dagger + U S \dd V^\dagger V V^\dagger \\
          &+ U_\perp \dd S_\perp V_\perp^\dagger + (\mathbb{1}_n - U U^\dagger) \dd U_\perp S_\perp V_\perp^\dagger + U_\perp S_\perp \dd V_\perp^\dagger (\mathbb{1}_m - V V^\dagger) \\
          &+ U S \dd V^\dagger (\mathbb{1}_m - V V^\dagger) + U U^\dagger \dd U_\perp S_\perp V_\perp^\dagger\\
          &+ (\mathbb{1}_n - U U^\dagger) \dd U S V^\dagger + U_\perp S_\perp \dd V_\perp^\dagger V V^\dagger.
  \end{align}
\end{subequations}
Applying the relations from Eqs.~(\ref{eq:diff_prop_2},~\ref{eq:diff_prop_3}) we obtain
\begin{subequations}
  \begin{align}
    \dd A &= U \dd S V^\dagger + U U^\dagger \dd U S V^\dagger + U S \dd V^\dagger V V^\dagger \\
          &+ U_\perp \dd S_\perp V_\perp^\dagger + (\mathbb{1}_n - U U^\dagger) \dd U_\perp S_\perp V_\perp^\dagger + U_\perp S_\perp \dd V_\perp^\dagger (\mathbb{1}_m - V V^\dagger) \\
          &+ U S \dd V^\dagger (\mathbb{1}_m - V V^\dagger) - U \dd U^\dagger U_\perp S_\perp V_\perp^\dagger\\
          &+ (\mathbb{1}_n - U U^\dagger) \dd U S V^\dagger - U_\perp S_\perp V_\perp^\dagger \dd V V^\dagger.
  \end{align}
\end{subequations}
To split up the equation into the different components, which has been proven to be practical, we multiply the equation from the left and right by the (semi)-unitaries and the projector into the discarded space,
\begin{align}
    U^\dagger \dd A V &= \dd S + U^\dagger \dd U S + S \dd V^\dagger V , \\
    (\mathbb{1}_n - U U^\dagger) \dd A (\mathbb{1}_m - V V^\dagger) &= U_\perp \dd S_\perp V_\perp^\dagger + (\mathbb{1}_n - U U^\dagger) \dd U_\perp S_\perp V_\perp^\dagger + U_\perp S_\perp \dd V_\perp^\dagger (\mathbb{1}_m - V V^\dagger) , \\
    U^\dagger \dd A (\mathbb{1}_m - V V^\dagger) &= S \dd V^\dagger (\mathbb{1}_m - V V^\dagger) - \dd U^\dagger U_\perp S_\perp V_\perp^\dagger ,  \label{eq:Ut_dA_1-VVt} \\
    (\mathbb{1}_n - U U^\dagger) \dd A V &= (\mathbb{1}_n - U U^\dagger) \dd U S - U_\perp S_\perp V_\perp^\dagger \dd V . \label{eq:1-UUt_dA_V}
\end{align}
Similar to the discussion in the previous Section~\ref{sec:svd_diff_square}, we define the following notation for compactness, which consists of the two terms generating $\dd U$ and $\dd V$ as stated in Eqs.~(\ref{eq:dU_full_gkl},~\ref{eq:dV_full_gkl}),
\begin{align}
    \widetilde{\dd U}_1 &= U^\dagger \dd U , & \widetilde{\dd U}_2 &= (\mathbb{1}_n - U U^\dagger) \dd U , \\
    \widetilde{\dd V}_1 &= V^\dagger \dd V , & \widetilde{\dd V}_2 &= (\mathbb{1}_m - V V^\dagger) \dd V.
\end{align}
It is easy to observe that the equations used to calculate $\dd S$, $\widetilde{\dd U}_1$ and $\widetilde{\dd V}_1$ are still the same so we can just use the result calculated in the sections before. Therefore, we just have to use the last two Eqs.~(\ref{eq:Ut_dA_1-VVt},~\ref{eq:1-UUt_dA_V}) to calculate the other two components. The first step is to observe that we can express the part $U_\perp S_\perp V_\perp^\dagger$ by multiplying Eq.~(\ref{eq:A_tSVD_gkl}) by $U_\perp U_\perp^\dagger = \mathbb{1}_n - UU^\dagger$ or $V_\perp V_\perp^\dagger = \mathbb{1}_m - VV^\dagger$, respectively,
\begin{align}
    U_\perp S_\perp V_\perp^\dagger &= (\mathbb{1}_n - U U^\dagger) A = A (\mathbb{1}_n - V V^\dagger) .
\end{align}
Inserting that into our two equations~(\ref{eq:Ut_dA_1-VVt},~\ref{eq:1-UUt_dA_V}), taking the Hermitian conjugate of Eq.~\eqref{eq:Ut_dA_1-VVt} and using our compact notation, we have the set of equations
\begin{align}
    (\mathbb{1}_m - V V^\dagger) \dd A^\dagger U &= \widetilde{\dd V}_2 S  - A^\dagger \widetilde{\dd U}_2, \\
    (\mathbb{1}_n - U U^\dagger) \dd A V &= \widetilde{\dd U}_2 S - A \widetilde{\dd V}_2.
\end{align}
Both equations can be solved for one of the terms,
\begin{align}
    \widetilde{\dd V}_2 S &= (\mathbb{1}_m - V V^\dagger) \dd A^\dagger U + A^\dagger \widetilde{\dd U}_2, \\
    \widetilde{\dd U}_2 S &= (\mathbb{1}_n - U U^\dagger) \dd A V + A \widetilde{\dd V}_2.
\end{align}
Assuming that $n \leq m$, it is more efficient to work with $A A^\dagger$ and thus calculating $\widetilde{\dd U}_2$ first. To this end, we multiply one equation by $S$ from the right and the other by $A A^\dagger$ from the left,
\begin{align}
    \widetilde{\dd U}_2 S^2 - A A^\dagger \widetilde{\dd U}_2 &= (\mathbb{1}_n - U U^\dagger) \dd A V S + A (\mathbb{1}_m - V V^\dagger) \dd A^\dagger U .
\end{align}
This equation is a Sylvester equation which can be solved numerically~\cite{10.1145/361573.361582,1102170}. The solution can then be used to calculate $\widetilde{\dd V}_2$ via the equation
\begin{align}
    \widetilde{\dd V}_2 &= (\mathbb{1}_m - V V^\dagger) \dd A^\dagger U S^{-1} + A^\dagger \widetilde{\dd U}_2 S^{-1} .
\end{align}

\subsubsection*{Summary of the results for iterative solver}

To summarize the results for the derivative of an iterative case, we have 
\begin{align}
    \dd S &= \mathbb{1}_t \circ \frac{1}{2} \left(U^\dagger \dd A V + V^\dagger \dd A^\dagger U\right) , \\
    \dd U &= U \widetilde{\dd U}_1 + \widetilde{\dd U}_2 , \\
    \dd V &= V \widetilde{\dd V}_1 + \widetilde{\dd V}_2 ,
\end{align}
with the terms
\begin{align}
    \widetilde{\dd U}_1 &= F \circ \left(U^\dagger \dd A V S + S V^\dagger \dd A^\dagger U\right) + S^{-1} \circ \frac{1}{2} \left(U^\dagger \dd A V - V^\dagger \dd A^\dagger U\right) , \\
    \widetilde{\dd V}_1 &= F \circ \left(S U^\dagger \dd A V + V^\dagger \dd A^\dagger U S\right) .
\end{align}
For $n \leq m$ the terms $\widetilde{\dd U}_2$ and $\widetilde{\dd V}_2$ can be found by solving a Sylvester equation,
\begin{align}
    \widetilde{\dd U}_2 S^2 - A A^\dagger \widetilde{\dd U}_2 &= (\mathbb{1}_n - U U^\dagger) \dd A V S + A (\mathbb{1}_m - V V^\dagger) \dd A^\dagger U , \\
    \widetilde{\dd V}_2 &= (\mathbb{1}_m - V V^\dagger) \dd A^\dagger U S^{-1} + A^\dagger \widetilde{\dd U}_2 S^{-1} .
\end{align}
For $n > m$ the equation would be
\begin{align}
    \widetilde{\dd V}_2 S^2 - A^\dagger A \widetilde{\dd V}_2 &= (\mathbb{1}_m - V V^\dagger) \dd A^\dagger U S + A^\dagger (\mathbb{1}_n - U U^\dagger) \dd A V, \\
    \widetilde{\dd U}_2 &= (\mathbb{1}_n - U U^\dagger) \dd A V S^{-1} + A \widetilde{\dd V}_2 S^{-1}.
\end{align}

\section{Trunctated eigenvalue decomposition (EVD)}

The full eigenvalue decomposition of a general square matrix $A\in\mathbb{C}^{n\times n}$ given in terms of a right eigenvector matrix $X \in \mathbb{C}^{n \times n}$ and its corresponding diagonal eigenvalue matrix $\Lambda \in \mathbb{C}^{n \times n}$ can be stated as
\begin{align}
    A X &= X \Lambda .
\end{align}
We shall denote by $Y \equiv X^{-1} \in \mathbb{C}^{n \times n}$ the inverse of the eigenvector matrix $X$ for brevity. However, we occasionally write $X^{-1}$ explicitly when we deem it necessary to highlight specific points.

For the truncation of the eigenvalue decomposition into $p$ kept and $n - p$ discarded values we split the eigenvalues and eigenvectors into two blocks $\lambda/x/y$ and $\lambda_\perp/x_\perp/y_\perp$ such that
\begin{align}
    X &= \left(\begin{array}{c|c}x & x_\perp\end{array}\right), & \Lambda &= \left(\begin{array}{c|c}\lambda & 0 \\ \hline 0 & \lambda_\perp\end{array}\right), & Y &= \left(\begin{array}{c}y \\ \hline y_\perp\end{array}\right), \label{eq:eig_xy_rel}
\end{align}
with $x \in \mathbb{C}^{n \times p}$, $\lambda \in \mathbb{C}^{p \times p}$ diagonal, $y \in \mathbb{C}^{p \times n}$, $x_\perp \in \mathbb{C}^{n \times (n - p)}$, $\lambda_\perp \in \mathbb{C}^{(n - p) \times (n - p)}$ diagonal and $y_\perp \in \mathbb{C}^{(n - p) \times n}$.
Furthermore, we can find from the relation $X Y = Y X = \mathbb{1}_n$ the following equations:
\begin{align}
    y x &= \mathbb{1}_{p}, & y_\perp x_\perp &= \mathbb{1}_{n-p}, & y_\perp x &= y x_\perp = 0, & x y + x_\perp y_\perp = \mathbb{1}_n , \label{eq:eig_trunc_rel}
\end{align}
and for the eigenvalue relation
\begin{align}
    A x &= x \lambda, &  A x_\perp &= x_\perp \lambda_\perp, & A &= x \lambda y + x_\perp \lambda_\perp y_\perp .
    \label{eq:eig_trunc_val_rel}
\end{align}

\subsection{Full derivative}

The derivative of the full eigenvalue relation is
\begin{align}
    \dd A X + A \dd X &= \dd X \Lambda + X \dd \Lambda.
\end{align}
Premultipling by $X^{-1}$ and substitute $X^{-1}A = \Lambda X^{-1}$ yields
\begin{align}
    X^{-1} \dd A X - \dd \Lambda &= X^{-1} \dd X \Lambda - \Lambda X^{-1} \dd X .
\end{align}

As discussed in Ref.~\cite{aa07_comput_eigen_eigen_deriv_gener}, there is a gauge freedom in the eigenvalue decomposition such that the decomposition can be transformed by a gauge fixing matrix $\Gamma$. If we assume that we already have a set of eigenvectors $\bar{X}$ (for example calculated by LAPACK), we can relate both sets of eigenvectors by $\Gamma$,
\begin{align}
    X &= \bar{X} \Gamma , & \Lambda \Gamma &= \Gamma \Lambda .
\end{align}
With this relation we can write the derivative in terms of the known eigenvectors $\bar{X}$,
\begin{align}
    \Gamma^{-1} \bar{X}^{-1} \dd A \bar{X} \Gamma - \dd \Lambda &= \Gamma^{-1} \bar{X}^{-1} \dd X \Lambda - \Lambda \Gamma^{-1} \bar{X}^{-1} \dd X, \\
    \Rightarrow \bar{X}^{-1} \dd A \bar{X} \Gamma - \Gamma \dd \Lambda &= \bar{X}^{-1} \dd X \Lambda - \Lambda \bar{X}^{-1} \dd X.
\end{align}
In the following we want derive the derivative of the truncated eigenvalue decomposition from Eq.~\eqref{eq:eig_trunc_val_rel} starting from the general relation.

\subsection{Derivative for truncated EVD with non-degenerate eigenvalues}

In this section we assume that there are no degenerate eigenvalue and therefore the gauge-fixing matrix $\Gamma$ to be a diagonal matrix.

Applying our notation of the split up parts for the derivative, we obtain the four equations
\begin{align}
    \bar{y} \dd A \bar{x} \gamma - \gamma \dd \lambda &= \bar{y} \dd x \lambda - \lambda \bar{y} \dd x , \label{eq:eig_dlambda_full}\\
    \bar{y}_\perp \dd A \bar{x} \gamma &= \bar{y}_\perp \dd x \lambda - \lambda_\perp \bar{y}_\perp \dd x , \label{eq:eig_dx_full}\\
    \bar{y} \dd A \bar{x}_\perp \gamma_\perp &= \bar{y} \dd x_\perp \lambda_\perp - \lambda \bar{y} \dd x_\perp , \\
    \bar{y}_\perp \dd A \bar{x}_\perp \gamma_\perp - \gamma_\perp \dd \lambda_\perp &= \bar{y}_\perp\dd x_\perp \lambda_\perp - \lambda_\perp \bar{y}_\perp \dd x_\perp .
\end{align}
Since the last two equations are only related to the derivatives in the part we truncate, we are only interested in the first two equations. From the first equation we can directly obtain the result for the derivative of the eigenvalues, since this one is restricted to having only diagonal non-zero values. When we look at the diagonal elements of the first equation we observe that the right-hand side cancels out, thus finding
\begin{align}
    \dd \lambda &= \mathbb{1}_p \circ (\bar{y} \dd A \bar{x}).
\end{align}
The task at hand now is to determine $\dd x$. As before for the truncated SVD, we split the problem up into two parts, $\dd x = \bar{x} \bar{y} \dd x + (\mathbb{1}_n - \bar{x}\bar{y}) \dd x = \bar{x} \widetilde{\dd x}_1 + \widetilde{\dd x}_2$. The first part can be determined from Eq.~\eqref{eq:eig_dlambda_full}, the second one from Eq.~\eqref{eq:eig_dx_full}.
For the first part $\widetilde{\dd x}_1 = \bar{y} \dd x$, we can obtain the off-diagonal elements from Eq.~\eqref{eq:eig_dlambda_full},
\begin{align}
    \bar{\mathbb{1}}_p \circ \widetilde{\dd x}_1 &= F \circ \left( \bar{y} \dd A \bar{x} \gamma \right),
\end{align}
with $\bar{\mathbb{1}}$ being the matrix consisting of ones everywhere except for the zero-valued diagonal and $F \in \mathbb{C}^{p \times p}$ a matrix with elements $F_{i,j} = 1/(\lambda_j - \lambda_i)$ for $i \neq j$ and zero on the diagonal.
The diagonal elements are dependent on the gauge-fixing we choose. We assume that we set an $m$-th element of the $k$-th eigenvector to be $1$ (with $m$ not necessary the same for all eigenvectors) such that $\gamma_k = \frac{1}{x_{mk}}$. For details we want to point to Ref.~\cite{aa07_comput_eigen_eigen_deriv_gener} where the reason how to select the values of $m$ is discussed. Summarized, there are several methods such as choosing the element with largest magnitude of $x_k$ (LAPACK default) or the largest product of the absolute values of $x$ and $y$ such that $|x_{mk}||y_{km}| = \max_i |x_{ik}||y_{ki}|$~\cite{aa07_comput_eigen_eigen_deriv_gener}. However, since we fix this value, we know that $\dd x_{mk} = 0$. This can be used to express the diagonal element $\widetilde{\dd x}_{1_{kk}}$ in terms of the off-diagonal ones,
\begin{align}
    \dd x_{mk} &= \sum_l \bar{x}_{ml} \widetilde{\dd x}_{1_{lk}} + \widetilde{\dd x}_{2_{mk}} = 0, \\
    \Rightarrow \widetilde{\dd x}_{1_{kk}} &= - \left(\frac{1}{\bar{x}_{mk}} \sum_{l \neq k} \bar{x}_{ml} \widetilde{\dd x}_{1_{lk}} + \widetilde{\dd x}_{2_{mk}} \right)
\end{align}

To determine the second part, we have to remove the dependency on the $\perp$ elements in Eq.~\eqref{eq:eig_dx_full}. To this end, we can use, after multiplication with $\bar{x}_\perp$ from the left side, the relations we found in Eqs.~\eqref{eq:eig_trunc_rel} and \eqref{eq:eig_trunc_val_rel},
\begin{align}
    (\mathbb{1}_n - \bar{x}\bar{y}) \dd A \bar{x} \gamma &= (\mathbb{1}_n - \bar{x}\bar{y}) \dd x \lambda - (A - \bar{x} \lambda \bar{y}) \dd x .
\end{align}
With the relation $A \bar{x} \bar{y} = \bar{x} \lambda \bar{y}$ we can obtain a Sylvester equation for $\widetilde{\dd x}_2 = (\mathbb{1}_n - \bar{x}\bar{y}) \dd x$,
\begin{align}
    (\mathbb{1}_n - \bar{x}\bar{y}) \dd A \bar{x} \gamma &= \widetilde{\dd x}_2 \lambda - A \widetilde{\dd x}_2 .
\end{align}
After solving this equation numerically, we have all components together for the full eigenvector derivative.

\subsection{Derivative for degenerate eigenvalues}

In Ref.~\cite{aa07_comput_eigen_eigen_deriv_gener} the derivative for the case of degenerate eigenvalues has been discussed as well, but this requires the knowledge of the second differential $\dd^2 A$ for a unique solution. Since this object is not generally available in numerical implementation of automatic differentiation frameworks, the discussion of this case is beyond scope of this work.

\begin{acknowledgments}
    J.N.\ thanks Philipp Schmoll, Roberto Losada, Erik Lennart Weerda and Frederik Wilde for proof-reading this note and their very helpful comments. This work has been funded by the Deutsche Forschungsgemeinschaft (DFG, German Research Foundation) under the project number 277101999 – CRC 183 (project B01).
\end{acknowledgments}

\bibliography{references}

\end{document}